\documentclass[twoside]{article}

\usepackage{lipsum} 

\usepackage[sc]{mathpazo} 
\usepackage[T1]{fontenc} 
\usepackage{microtype} 

\usepackage{amsmath,graphicx,amssymb,amsfonts}

\usepackage{graphicx}
\usepackage{tkz-berge}
\usepackage{tikz}
\usepackage{tkz-graph}
\usetikzlibrary{fit,positioning,calc}
\usetikzlibrary{shapes,arrows}

\graphicspath{{images/}{../images/}}

\usepackage[hmarginratio=1:1,top=32mm,columnsep=20pt]{geometry} 
\usepackage{multicol} 
\usepackage[hang, small,labelfont=bf,up,textfont=it,up]{caption} 
\usepackage{booktabs} 
\usepackage{float} 
\usepackage{hyperref} 

\usepackage{lettrine} 
\usepackage{paralist} 

\usepackage{abstract} 

\usepackage{titlesec} 
\renewcommand\thesection{\arabic{section}} 
\renewcommand\thesubsection{\thesection.\arabic{subsection}} 
\titleformat{\section}[block]{\large\scshape\centering}{\thesection.}{1em}{} 
\titleformat{\subsection}[block]{\large}{\thesubsection.}{1em}{} 

\usepackage{fancyhdr} 
\newtheorem{thm}{Theorem}[subsection]
\newtheorem{cor}[thm]{Corollary}

\newtheorem{lem}[thm]{Lemma}

\title{\vspace{-15mm}\fontsize{24pt}{10pt}\selectfont\textbf{ Domination in commuting graph and its complement }} 

\author{
\large
\textsc{E. Vatandoost and Y. Golkhandypour}\\[2mm] 
\normalsize Department of mathematics, International university of imam Khomeini, Qazvin, Iran \\ 
\normalsize \href{vatandoost@sci.ikiu.ac.ir and y.golkhandypour@edu.ikiu.ac.ir}{vatandoost@sci.ikiu.ac.ir and y.golkhandypour@edu.ikiu.ac.ir} 
\vspace{-5mm}
}
\date{}

\begin{document}

\maketitle 

\thispagestyle{fancy} 


\begin{abstract}
For each non-commutative ring R, the commuting graph of R is a graph with vertex set $R\setminus Z(R)$ and two vertices $x$ and $y$ are adjacent if and only if $x\neq y$ and $xy=yx$. In this paper, we consider the domination and signed domination numbers on commuting graph $\Gamma(R)$ for non-commutative ring $R$ with $Z(R)=\{0\}$. For a finite ring $R$, it is shown that $\gamma(\Gamma(R)) + \gamma(\overline{\Gamma}(R))=|R|$ if and only if $R$ is non-commutative ring on 4 elements. Also we determine the domination number of $\Gamma(\prod_{i=1}^{t}R_i)$ and commuting graph of non-commutative ring $R$ of order $p^3$, where $p$ is prime. Moreover we present an upper bound for signed domination number of $\Gamma(\prod_{i=1}^{t}R_i)$.\\
{\bf keywords:} Noncommutative ring; Commuting graph; Domination number; signed domination number.
\end{abstract}


\begin{multicols}{2} 

\section{Introduction}

\lettrine[nindent=0em,lines=3]{L} et $R$ be a non-commutative ring, $Z(R)$ denoted the center of $R$ and for $a\in R$, $C(a)$ denotes the centeralizer of $a$ in $R$. The commuting graph of $R$, denoted by $\Gamma(R)$, is a graph with vertex set $R\setminus Z(R)$ and joined two vertices $x$ and $y$ if and only if $x\neq y$ and $xy=yx$. This graph was introduced by Akbari et al. \cite{122}, and the complement of commuting graph of $R$ is denoted by $\overline{\Gamma}(R)$. See \cite{144}, \cite{122}, \cite{119}, \cite{145} and \cite{A1} for more details.\\
Let $G = (V,E)$ be a simple graph comprising a set $V(G)$ of vertices together with a set $E(G)$ of edges. A graph $G$ is said to be $connected$ if each pair of vertices are joined by a walk. The number of edges of the shortest walk joining $v$ and $u$ is called the $distance$ between $v$ and $u$ and denoted by $d(v,u)$. The maximum value of the distance function in a connected graph $G$ is called the $diameter$ of $G$ and denoted by $diam(G)$. The neigbours of a vertex $v\in V(G)$ is the set of edges incident to the $v$ and denoted by $N(v)$ and $|N(v)|=deg(v)$. The maximum degree of a graph $G$ denoted by $\Delta(G)$, and the minimum degree of a graph $G$ denoted by $\delta(G)$, are the maximum and minimum degree of its vertices. The $complete~graph$, $K_n$, is a graph with $n$ vertices in which each pair of vertices are adjacent. A {\it complete bipartite graph} is a graph whose vertices can be divided into two disjoint set $A$ and $B$ such that each edge is incident to a vertex in $A$ and a vertex in $B$ and denoted by $K_{n,m}$ where $|A|=n$ and $|B|=m$.\\
The {\it strong product}, $G\boxtimes H$, of graphs $G$ and $H $ is a graph whose structure is as follows:\\
{\bf i)} The vertex set of $G\boxtimes H$ is the Cartesian product $V(G)\times V(H)$.\\
{\bf ii)} Any two distinct vertices $(v,v')$ and $(u,u')$ are adjacent in $G\boxtimes H$ if and only if $v$ is adjacent to $u$ in $G$ and $v'=u'$, or $v=u$ and $v'$ is adjacent to $u'$ in $H$, or $v$ is adjacent to $u$ in $G$ and $v'$ is adjacent to $u'$ in $H$.\\
The {\it corona} $G=G_1\circ G_2$ is the graph formed from one copy of $G_1$ and $|V(G_1)|$ copies of $G_2$ where the $i$th vertex of $G_1$ is adjacent to every vertex in the $i$th copy of $G_2$. \\
A {\it dominating set} of G is a subset D of V(G) such that every vertex not in D is adjacent to at least one vertex in D.
The domination number of $G$ is the number of vertices in a minimal dominating set for G and denoted by $\gamma(G)$. See \cite{111}, \cite{121} and \cite{113} for more details.\\
The closed neighbour of $v$, denoted by $N[v]$, is the set $N(v)\cup \{v\}$. A function $f: V\rightarrow \{-1,1\}$ is a signed dominating function if for every vertex $v\in V(G)$, the closed neighbour of $v$ contains more vertices with function value 1 than with -1. Thus $f$ is a signed dominating function if $f[v]\geq 1$ for all $ v\in V(G)$, where $f[v] = \sum_{v\in V(G)} f(v)$. The weight
of $f$, denoted $f(G)$, is the sum of the function value of all vertices in $G$,
i.e., $f(G) = \sum_{x\in V(G)} f(x)$ . The {\it signed domination number} of $G$, denoted $\gamma_s(G)$, is the
minimum weight of signed dominating functions of $G$. Also the set of vertices with function value $-1$ is denoted by $V^-(G)$.\\
In this paper, we consider the domination and signed domination numbers on commuting graph $\Gamma(R)$ for non-commutative ring $R$ with $Z(R)=\{0\}$. For a finite ring $R$, it is shown that $\gamma(\Gamma(R)) + \gamma(\overline{\Gamma}(R))=|R|$ if and only if $R$ is the non-commutative ring on 4 elements. Also we determine the domination number of $\Gamma(\prod_{i=1}^{t}R_i)$ and commuting graph of non-commutative ring $R$ of order $p^3$, where $p$ is prime. Moreover we present an upper bound for the signed domination number of $\Gamma(\prod_{i=1}^{t}R_i)$. The main results in this paper are the following.\\
\\
{\bf Theorem A.} Let $R$ be a non-commutative ring of order $n$ and $Z(R)=\{0\}$. Then\\
{\it i)} $\gamma(\Gamma(R)) + \gamma(\overline{\Gamma}(R))=n$ if and only if $R$ is isomorphic with one of the following rings:\\
$E=\langle x,y \: : \: 2x=2y=0 \: , \: x^2=x \: , \:y^2=y \: , \: xy=x , \: yx=y \rangle$\\
$F=\langle x,y \: : \: 2x=2y=0 \: , \: x^2=x \: , \: y^2=y \: , \: xy=y , \: yx=x \rangle$.\\
{\it ii)} $\gamma(\Gamma(R)) + \gamma(\overline{\Gamma}(R))\neq n-1$.\\
{\it iii)} $\gamma(\Gamma(R)) + \gamma(\overline{\Gamma}(R))=n-2$ if and only if $n$ is even and $\Gamma(R)=K_3\cup (n-4)K_1$.\\
\\
{\bf Theorem B.} Let R be a non-commutative ring of order $p^3$ and $Z(R)=\{0\}$. Then\\
{\it i)} $\gamma(\Gamma(R))=p^2+p+1$.\\
or\\
{\it ii)} $\gamma(\Gamma(R))=\ell_1+\ell_2$, where $\ell_1$ and $\ell_2$ satisfy in $\ell_1+(p+1)\ell_2=p^2+p+1$.\\
\\
{\bf Theorem C.} Let $R_i$, $(1\leq i\leq t)$, be non-commutative ring of order $n_i$ and $Z(R_i)=\{0\}$. Then
$\gamma(\Gamma(\prod_{i=1}^{t}R_i))=Min_{1\leq i\leq t}(\gamma(\Gamma(R_i)))$. \\
\\
{\bf Theorem D.} Let $R$ be a non-commutative ring of order $n$ and $Z(R)=\{0\}$.\\
{\it i)} Let n be even. Then $\gamma_s(\Gamma(R))=n-1$ if and only if $R$ is isomorphic with one of the following rings:\\
$E=\langle x,y~|~2x=2y=0, x^2=x, y^2=y, xy=x, yx=y\rangle$\\
$F=\langle x,y~|~2x=2y=0, x^2=x, y^2=y, xy=y, yx=x\rangle$.\\
 {\it ii)} Let n be odd. Then $\gamma_s(\Gamma(R))=n-1$ if and only if $\Gamma(R)$ is the union of $\frac{n-1}{2}$ copies of $P_2$.\\
 \\
 {\bf Theorem E.} Let $R_i$, $1\leq i\leq t$ be non-commutative ring such that $|R_i|=n_i$ and $Z(R_i)=\{0\}$. Also, let $\delta_i$ be minimum degree of $\Gamma(R_i)$.\\
{\it i)} If for all $1\leq i\leq t$, $\delta_i$ is odd, then $\gamma_s(\Gamma(\prod_{i=1}^{t}R_i))\leq \prod_{i=1}^{t}n_i-\prod_{i=1}^t(\delta_i+2)+2$.\\
{\it ii)} Otherwise, $\gamma_s(\Gamma(\prod_{i=1}^{t}R_i))\leq \prod_{i=1}^{t}n_i-\prod_{i=1}^t(\delta_i+2)+1$.


\section{\bf Preliminaries}
First we give some facts that are needed in the section III.

\subsection{ On commuting graph}
\begin{lem}\label{3} Let $R$ be a non-commutative ring of order $n$ and $Z(R)=\{0\}$. If $\overline{\Gamma}(R)$ contains a vertex of degree $k$, then $k>\frac{n-1}{2}$.
\end{lem}
{\it Proof.}
On the contrary, let $v$ be a vertex of degree $k$ in $\overline{\Gamma}(R)$ such that $k\leq \frac{n-1}{2}$. So $|C(v)|\geq \frac{n-1}{2}+1$. Thus $|C(v)| \nmid n$, which is impossible. Hence $deg(v)>\frac{n-1}{2}$, for each $v\in \overline{\Gamma}(R)$.
$\hfill\Box$
\begin{lem}\label{30} Let $R$ be a non-commutative ring of order $n$ and $Z(R)=\{0\}$. Then $\Gamma(R)$ is not a cycle. Also, $\Gamma(R)$ does not have $C_4$ as a component.
\end{lem}
{\it Proof.} On the contrary, let $\Gamma(R)=C_{n-1}$ and $v_1, \ldots, v_{n-1}\in V(C_{n-1})$, such that $v_i\in N(v_{i+1})$. So $C(v_i)=\{0, v_{i-1}, v_i, v_{i+1}\}$. Since $C(v_i)$ is a subgroup of $(R,+)$, $v_{i-1}+v_{i+1}\in C(v_i)$. If $v_{i-1}+v_{i+1}=0$, then $v_{i-1}=-v_{i+1}$. It follows that $v_{i-2}\in N(v_{i+1})$, which is impossible. If $v_{i-1}+v_{i+1}\in \{v_{i-1}, v_{i+1}\}$, then $v_{i-1}=0$ or $v_{i+1}=0$, a contradiction. Thus $v_{i-1}+v_{i+1}=v_i$. Also, $C(v_{i-1})=\{0, v_{i-1}, v_{i-2}, v_{i}\}$. Similarly, $v_{i-2}+v_{i-1}=v_{i}$. Hence $v_{i-2}=v_{i+1}$, which is a contradiction. Therefore $\Gamma(R)$ is not a cycle. \\
Also let $\Gamma(R)$ has $C_4$ as a component and let $V(C_4)=\{x_1, x_2, x_3, x_4\}$ such that $x_1\notin N(x_2)$. Then $|C(x_1)\cap C(x_2)|=3$, which is impossible.
$\hfill\Box$
\begin{lem}\label{40} Let $R$ be a non-commutative ring of order $n$ and $Z(R)=\{0\}$. Then $\Gamma(R)$ does not have both an isolated vertex and a vertex of degree one.
\end{lem}
{\it Proof.} On the contrary, let $x$ be an isolated vertex and $y$ be a vertex of degree 1 in $\Gamma(R)$. So $C(x)=\{0, x\}$ and $C(y)=\{0, y, -y\}$. If $C(x)+C(y)=H$, then $H$ is a subgroup of $(R,+)$ and $|H|=6$. Since $x+y\in H$ and $O(x+y)=6$, $(H,+)$ is a cyclic group. Thus $C(x+y)=H$. It follows that $K_5$ is a subgraph of commuting graph. So $x$ is not an isolated vertex, which is contradiction.
$\hfill\Box$
\begin{thm}\cite{A1}\label{20} Let $R$ be a non-commutative finite ring with $|R|>4$. Then $diam(\overline{\Gamma}(R))=2$.
\end{thm}
\begin{thm}\cite{A1}\label{21} Let $R$ be a non-commutative ring. Then $\overline{\Gamma}(R)$ is not a complete bipartite graph.
\end{thm}
\begin{lem}\label{16} Let $R$ be a non-commutative ring of order $n$ and $Z(R)=\{0\}$. Then $n\neq 6$.
\end{lem}
{\it Proof.}
On the contrary, let $n=6$. By Lemma \ref{3}, $\delta(\overline{\Gamma}(R))\geq 3$. If $\delta(\overline{\Gamma}(R))=4$, then $\overline{\Gamma}(R)=K_5$, contrary to Theorem \ref{20}. So $\delta(\overline{\Gamma}(R))=3$. Since there is no $3-$regular graph on 5 vertices, $\Delta(\overline{\Gamma}(R))=4$. Hence $\overline{\Gamma}(R)$ is one of the graphs that are shown in figure 1.
  \begin{center}
    \includegraphics[width=0.48\textwidth]{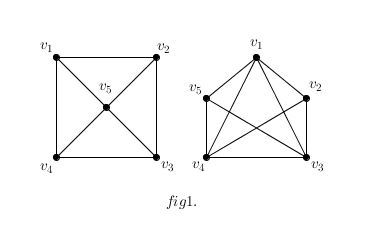}
  \end{center}
In both graphs, $\Gamma(R)$ is a union of isolated vertices and copies of $P_2$, contrary to Lemma \ref{40}. Therefore $n\neq 6$.
$\hfill\Box$
\begin{cor}\label{4} Let $R$ be a non-commutative ring of order $n$ and $Z(R)=\{0\}$. Then $n=4$ or $n\geq 8$.
\end{cor}
\begin{lem}\label{6} Let R be a ring of order $p$, where p is prime. Then R is a commutative ring.
\end{lem}
{\it Proof.}
The proof is straightforward.
$\hfill\Box$
\begin{thm}\label{2}\cite{A1} Let R be a non-commutative ring of order $p^2$, where $p$ is a prime number. Then R is one of the following rings. \\
$E=\langle x,y \: : \: px=py=0 \: , \: x^2=x \: , \: y^2=y \: , \: xy=x , \: yx=y \rangle$\\
$F=\langle x,y \: : \: px=py=0 \: , \: x^2=x \: , \: y^2=y \: , \: xy=y , \: yx=x \rangle$.
\end{thm}
\begin{thm}\cite{A1}\label{18} Let $R$ be a non-commutative finite ring with $diam(\overline{\Gamma}(R))=1$. Then $R$ is of type $E$ or $F$ (see Theorem \ref{2}).
\end{thm}
\begin{thm}\cite{122}\label{17} For any non-commutative ring $R$ and $x,y\in V(\overline{\Gamma}(R))$, there is a path between $x$ and $y$ in $\overline{\Gamma}(R)$ whose length is at most two.
\end{thm}
\begin{lem}\label{10} Let $R$ be a finite ring of order $p^2$ and $Z(R)\neq \{0\}$. Then $R$ is commutative.
\end{lem}
{\it Proof.}
On the contrary, let $R$ be a non-commutative ring. It follows immediately that $|Z(R)|=p$. Clearly, for any $x\in R\setminus Z(R)$, $|C(x)|>p$ and $|C(x)|\mid p^2$. So $|C(x)|=p^2$ and $x\in Z(R)$, a contradiction. Therefore $R$ is commutative ring.
$\hfill\Box$
\begin{lem}\label{8} Let $R$ be a non-commutative ring and $x,y\in V(\Gamma(R))$ such that $C(x)$ and $C(y)$ are commutative. If $y\in C(x)$, then $C(x)=C(y)$.
\end{lem}
{\it Proof.} The proof is straightforward.
$\hfill\Box$
\begin{thm} \label{44}\cite{146} {\bf (Scorza)} Let $\{A_i: 1\leq i\leq 3\}$ be an irredundant cover with core-free intersection $D$ for a group $G$. Then $D=1$ and $G\cong Z_2\oplus Z_2$.
\end{thm}
\subsection{On domination number}
\begin{lem}\label{19}\cite{123} Let $G$ be a graph on $n$ vertices. Then $\gamma(G)=1$ if and only if $\Delta(G)=n-1$.
\end{lem}
\begin{thm}\cite{124}\label{Gbar} Let $G$ be a graph on $n$ vertices. Then
\item[i)] $\gamma(G)+\gamma(\overline{G}) \leq n+1$.
\item[ii)] $\gamma(G)\gamma(\overline{G}) \leq n$.
\end{thm}
\begin{thm}\cite{125}\label{c4} Let $G$ be a graph without isolated vertices on $n$ vertices such that $n$ is even. Then $\gamma(G)= \frac{n}{2}$ if and only if the components of $G$ are $C_4$ or $H\circ K_1$ where $H$ is a connected graph.
\end{thm}
\begin{thm} \cite{113} \label{39} Let $G$ be a graph with no isolated vertex. Then $\gamma(G)\leq \frac{n}{2}$.
\end{thm}
\begin{thm}\label{13}\cite{116}, \cite{114} For any graph $G$,
\begin{center}
 $\lceil \frac{n}{1+\Delta(G)} \rceil\leq \gamma(G)\leq n-\Delta(G)$.
\end{center}
\end{thm}
\begin{thm}\label{14}\cite{115} If a graph $G$ has no isolated vertices, then
\begin{center}
$\gamma(G)\leq \frac{n+2-\delta(G)}{2}$.
\end{center}
\end{thm}
\subsection{On signed domination number}
\begin{lem}\label{1}\cite{118} A graph $G$ has $\gamma_s (G)=n$ if and only if every $v\in G$ is either isolated, an endvertex or adjacent to an
endvertex.
\end{lem}
\begin{lem} Let $G$ be a graph on $n$ vertices and $\alpha$ be an odd number. Then $\gamma_s(G)\neq n-\alpha$.
\end{lem}
{\it Proof.} The proof is straightforward.
$\hfill\Box$
\begin{thm}\cite{117} \label{9} Let $K_n$ be a complete graph on $n$ vertices. Then
\begin{equation*}
\gamma_s(K_n)=
      \begin{cases}
       2 & n ~be ~even\\
        1& n~be~ odd.
      \end{cases}
\end{equation*}
\end{thm}
\begin{lem}\label{5} Let $G$ be a graph with $\delta(G)\geq 6$. Then $|V^{-}(G)|\geq 3$.
\end{lem}
{\it Proof.}
Let $v\in V(G)$ and $deg(v)=\delta(G)\geq 6$. If $N(v)=\{v_1, \ldots, v_{\delta(G)}\}$, then consider the function $f: V(G)\rightarrow \{-1, 1\}$ for which $f(v)=f(v_1)=f(v_2)=-1$. Clearly, $f[w]\geq 1$ where $w\in \{v, v_1, v_2\}$. So $|V^{-}(G)|\geq 3$.
$\hfill\Box$
\begin{lem}\label{26} Let $R$ be a non-commutative ring of order 8 and $Z(R)=\{0\}$. Then $\gamma_s(\overline{\Gamma}(R))=1$.
\end{lem}
{\it Proof.}
Let $v\in V(\overline{\Gamma}(R))$ such that $deg(v)=k$. By Lemma \ref{3}, $k\geq 4$. If $k=5$, then $|C(v)|\nmid |R|$, which is a contradiction. Hence $k\in \{4, 6\}$. If $\delta(\overline{\Gamma}(R))=6$, then $\overline{\Gamma}(R)=K_7$, contrary to Theorem \ref{20}. \\
Let $\delta(\overline{\Gamma}(R))=4$. Then $\Delta(\overline{\Gamma}(R))\in \{4, 6\}$. If $\Delta(\overline{\Gamma}(R))=4$, then $\overline{\Gamma}(R)$ is a 4-regular graph on 7 vertices, which are depicted in figure 2.
  \begin{center}
    \includegraphics[width=0.48\textwidth]{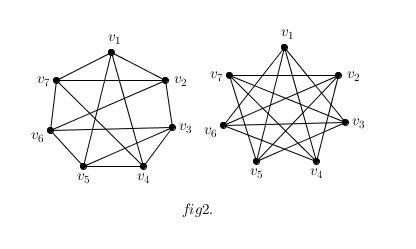}
  \end{center}
In both graphs, $\Gamma(R)=C_7$, contrary to Lemma \ref{30}. \\
Let $\Delta(\overline{\Gamma}(R))=6$. Then $\overline{\Gamma}(R)$ has at least two vertices $v$ and $u$ of degree 4. We consider the following two cases.\\
{\bf Case 1.}  Let $v\in N(u)$ in $\overline{\Gamma}(R)$. Then $C(v)=\{0, v, x_1, x_2\}$ and $C(u)=\{0, u, y_1, y_2\}$. So $|C(v)\cap C(u)|\in \{1, 2\}$.\\
If $|C(v)\cap C(u)|=1$, then $|C(v)+C(u)|>|R|$, which is impossible. \\
Let $|C(v)\cap C(u)|=2$, $C(v)=\{0, v, x_1, x_2\}$ and $C(u)=\{0, u, x_1, y_2\}$. Thus $x_1\notin N(v)\cup N(u)$, $y_2\in N(v)$ and $x_2\in N(u)$ in $\overline{\Gamma}(R)$ and there exist $z_1, z_2\in V(\overline{\Gamma}(R))$ such that $z_1, z_2\in N(v)\cap N(u)$. This case will be divided into 2 subcases. \\
{\bf Subcase i.} If $deg(x_1)=6$, then $x_1\in N(v)\cap N(u)$, a contradiction. \\
{\bf Subcase ii.} Let $deg(x_1)=4$. Since $x_1\notin N(v)\cup N(u)$, it follows that $\{y_2, x_2, z_1, z_2\}\subseteq N(x_1)$. Similarly, since $y_2\in N(v)$ and $x_2\in N(u)$, it follows that $\{z_1, z_2\}\subseteq N(y_2)\cap N(x_2)$. Thus $deg(x_2)=deg(y_2)=4$. We claim that $z_1\notin N(z_2)$ in $\overline{\Gamma}(R)$. For let $z_1\in N(z_2)$. Then by the above argument, $z_1$ and $z_2$ have exactly two common neigbours. But $v, u, y_2, x_2\in N(z_1)\cap N(z_2)$, which is a contradiction. Hence $z_1\notin N(z_2)$, as claimed. Thus $deg(z_1)=deg(z_2)=5$, which is impossible. Therefore this case will not happen. \\
{\bf Case 2.} Let $v\notin N(u)$ in $\overline{\Gamma}(R)$. Then $C(v)=C(u)=\{0, v, u, \alpha\}$. So there exist $z_1, z_2, z_3, z_4\in V(\overline{\Gamma}(R))$ such that $z_1, z_2, z_3, z_4\in N(v)\cap N(u)$ in $\overline{\Gamma}(R)$. \\
If $deg(\alpha)=6$ in $\overline{\Gamma}(R)$, then $\alpha\in N(v)$, a contradiction. So $deg(\alpha)=4$ and $z_1, z_2,  z_3, z_4\in N(\alpha)$. We claim that $deg(z_i)=6$, where $1\leq i\leq 4$. On the contrary, let there exist $1\leq i\leq 4$ such that $deg(z_i)=4$. Since $z_i\in N(v)$, as in the Case 1, $z_i$ and $v$ should have two common neigbours. Hence $deg(z_i)=5$, which is a contradiction. Thus $deg(z_i)=6$. \\
Consider the function $f: V(\overline{\Gamma}(R))\rightarrow \{-1, 1\}$ for which $f(v)=f(u)=f(\alpha)=-1$ and $f(z_i)=1$, where $1\leq i\leq 4$. Clearly, $f[w]\geq 1$, for all $w\in V(\overline{\Gamma}(R))$. This implies that $|V^{-}(\overline{\Gamma}(R))|\geq 3$.\\
If $g$ is a signed domination function on $\overline{\Gamma}(R)$ and $|\{w; g(w)=-1\}|>3$, then $g[\ell]=-1$ for some $\ell\in \{z_1, z_2, z_3, z_4\}$, which is a contradiction. Therefore $|V^-(\overline{\Gamma}(R))|=3$, and this completed the proof.
$\hfill\Box$
\begin{lem}\label{24} Let $R$ be a non-commutative ring of order $2p$, where $p$ is an odd prime and $Z(R)=\{0\}$. Then $\gamma_s(\overline{\Gamma}(R))=2$.
\end{lem}
{\it proof.} Let $|R|=2p$. Then by Lemma \ref{16}, $p>3$. Let $v\in V(\overline{\Gamma}(R))$ such that $deg(v)=k$. By Lemma \ref{3}, $k\geq p$. Clearly, if $k\notin \{p, 2p-2\}$, then $|C(v)|\nmid 2p$, which is impossible. Thus $k\in \{p, 2p-2\}$. If $\delta(\overline{\Gamma}(R))=2p-2$, then $\overline{\Gamma}(R)$ is a complete graph, contrary to Theorem \ref{20}. If $\Delta(\overline{\Gamma}(R))=p$, then we have a $p$-regular graph on $2p-1$ vertices, which is impossible. So $\Delta(\overline{\Gamma}(R))=2p-2$ and $\overline{\Gamma}(R)$ has at least two vertices $v$ and $u$ of degree $p$. The following two cases will be considered.\\
{\bf Case 1.}  Let $v\in N(u)$ in $\overline{\Gamma}(R)$. Then $C(v)=\{0, v, x_1, \ldots, x_{p-2}\}$ and $C(u)=\{0,u,y_1,\ldots,y_{p-2}\}$. Thus
$|C(v)\cap C(u)|=1$ and so $|C(v)+C(u)|=p^2>2p$, which is impossible. \\
Therefore this case will not happen.\\
{\bf Case 2.} Let $v\notin N(u)$ in $\overline{\Gamma}(R)$. Then $C(v)=\{0, v, u, x_1, \ldots, x_{p-3}\}$ and $C(u)=\{0, u, v, y_1, \ldots, y_{p-3}\}$. Thus $|C(v)\cap C(u)|=p$ and so  $C(v)=C(u)=\{0, v, u, z_1, z_2, \ldots, z_{p-3}\}$. Hence there exist $\alpha_1, \ldots, \alpha_p\in V(\overline{\Gamma}(R))$ such that $\alpha_1, \ldots, \alpha_p\in N(v)\cap N(u)$ in $\overline{\Gamma}(R)$. Obviously, $z_1, z_2, \ldots, z_{p-3}\notin N(v)\cup N(u)$. We claim that induced subgraph on $\{z_1, z_2, \ldots, z_{p-3}\}$ is empty. On the contrary, let $z_1\in N(z_2)$. By Case 1, $deg(z_1)\neq p$, $deg(z_2)\neq p$ or $deg(z_1),deg(z_2)\neq p$. Whithout loss of generality, let $deg(z_1)=2p-2$. Since $|\{z_2, z_3, \ldots, z_{p-3}\}|+|\{\alpha_1, \ldots, \alpha_p\}|=2p-4$, it follows that $deg(z_1)\leq 2p-4$, which is a contradiction. Thus $z_1\notin N(z_2)$. Similarly, $z_i\notin N(z_j)$, where $1\leq i,j \leq p-3$, and so induced subgraph on $\{z_1, z_2, \ldots, z_{p-3}\}$ is empty, as claimed. Hence $deg(z_i)=p$, where $1\leq i\leq p-3$ and by Case 1, for each $\ell\in \{\alpha_1, \ldots, \alpha_p\}$, $deg(\ell)=2p-2$. \\
Consider the function $f: V(\overline{\Gamma}(R))\rightarrow \{-1, 1\}$ for which $f(v)=f(u)=f(z_i)=-1$, where $1\leq i\leq p-3$ and $f(\alpha_j)=1$, where $1\leq j\leq p$. Clearly, $f[w]\geq 1$, for all $w\in V(\overline{\Gamma}(R))$. This implies that $|V^{-}(\overline{\Gamma}(R))|\geq p-1$. \\
If $g$ is a signed domination function on $\overline{\Gamma}(R)$ and $|\{w; g(w)=-1\}|>p-1$, then $g[\ell]=-1$ for some $\ell\in \{\alpha_1, \ldots, \alpha_p\}$, which is contradiction. Therefore $|V^{-}(\overline{\Gamma}(R))|=p-1$ and this completed the proof.
$\hfill\Box$
\section{\bf Main Results}
In this section we prove our main results.
\subsection{Domination number in $\Gamma(R)$}
\begin{thm}\label{3k} Let $R$ be a non-commutative ring of order odd number $n$ and $Z(R)=\{0\}$. Then $\gamma(\Gamma(R))=\frac{n-1}{2}$ if and only if $\Gamma(R)$ is a union of $\frac{n-1}{2}$ copies of $P_2$. Also $n=3^k$ for $k >1$.
\end{thm}
{\it Proof.} Let $\gamma(\Gamma(R))= \frac{n-1}{2}$. Since $n$ is odd, $\Gamma(R)$ does not have isolated vertex. By Theorem \ref{c4}, the components of $\Gamma(R)$ are $C_4$ or $H\circ K_1$ where $H$ is a connected graph. By Lemma \ref{30}, $\Gamma(R)$ does not have $C_4$ as a component. Hence the components of $\Gamma(R)$ are $H \circ K_1$. Let $x$ be an endvertex in $\Gamma(R)$ and $x\in N(y)$. Then $C(x)=\lbrace 0,x,y \rbrace$. If $O(x)=2$, then $n$ is even, which is false. So $O(x)\neq 2$. If $-x\neq y$, then $x$ is adjacent to $-x$ in $\Gamma(R)$, which is a contradiction. Hence $-x=y$ and so $deg(y)=1$. Thus $\Gamma(R)$ is union of $\frac{n-1}{2}$ copies of $P_2$. Therefore for every $0 \neq z \in R$, $C(z)=\lbrace 0,z,-z \rbrace$ and so $O(z)=3$. Hence $n=3^k$ for $k>1$. The proof of converse is easy.
$\hfill\Box$
\begin{cor} Let $R$ be a non-commutative ring of order $n$ such that $n$ and $|Z(R)|=t$ are odd. Then $\gamma(\Gamma(R))=\frac{n-t}{2}$ if and only if $\Gamma(R)$ is the union of $\frac{n-t}{2}$ copies of $P_2$.
\end{cor}
\begin{lem}\label{37} Let $R$ be a non-commutative ring of order $n$ and $Z(R)=\{0\}$. Then $\gamma(\Gamma(R)) \geq 3$.
\end{lem}
{\it Proof.} On the contrary, let $\gamma(\Gamma(R)) < 3$. If $\gamma(\Gamma(R))=1$, then there exist $x \in R \setminus Z(R)$ such that $C(x)=R$, which is impossible. Let $\gamma(\Gamma(R))=2$ and $D=\lbrace x,y \rbrace$ be a dominating set in $\Gamma(R)$. Then $R=C(x) \cup C(y)$. Hence $C(x) \subseteq C(y)$ or $C(y) \subseteq C(x)$. Without loss of generality, let $C(x) \subseteq C(y)$. Then $R=C(y)$ and so $y \in Z(R)$, which is a contradiction.
$\hfill\Box$
\begin{lem}\label{38} Let $R$ be a non-commutative ring of order $n$ and $Z(R)=\{0\}$. Then $\gamma(\Gamma(R)) =3$ if and only if $R$ is isomorphic with one of the following rings:\\
$E=\langle x,y \: : \: 2x=2y=0 \: , \: x^2=x \: , \: y^2=y \: , \: xy=x , \: yx=y \rangle$\\
$F=\langle x,y \: : \: 2x=2y=0 \: , \: x^2=x \: , \: y^2=y \: , \: xy=y , \: yx=x \rangle$.
\end{lem}
{\it Proof.} It is not hard to see that $\gamma(\Gamma(E))=\gamma(\Gamma(F))=3$. Conversely, let $\gamma(\Gamma(R))=3$ and $D=\{x,y,z \}$ be a dominating set in $\Gamma(R)$. Then $R=C(x) \cup C(y) \cup C(z)$. By Theorem \ref{44}, $(R,+)\cong Z_2 \oplus Z_2$. On the other hand, exactly $E$ and $F$ are non-commutative rings between all rings of order four.
$\hfill\Box$
\begin{cor} Let $R$ be a non-commutative ring of order $n$ and $Z(R)\neq \{0\}$. Then $\gamma(\Gamma(R))\geq 4$.
\end{cor}
\begin{lem} Let $R$ be a non-commutative ring of order $n$ with $Z(R)=\{0\}$. If $\gamma(\overline{\Gamma}(R))=1$, then $n=2^t$ for positive integer $t$.
\end{lem}
{\it Proof.} Let $D=\{x\}$ be a dominating set in $\overline{\Gamma}(R)$. So $x$ is an isolated vertex in $\Gamma(R)$, and so $O(x)=2$. Hence $n=2k$. On the contrary, let $p\mid n$, where $p$ is an odd prime. So there exist $y\in R$ such that $O(y)=p$. Hence $py=0,px=x$ and $2y\neq 0$. So $2ypx=0$. Thus $2y \in C(x)$ and so $2y=x$. It follows that $O(y)=4$, which is a contradiction. Therefore $n=2^t$.
$\hfill\Box$
\begin{lem}\label{42} Let $R$ be a non-commutative ring of order odd number $n$, $Z(R)=\{0\}$ and $3\nmid n$. Then $3<\gamma(\Gamma(R))<\frac{n-1}{2}$.
\end{lem}
{\it Proof.}
Since $n$ is odd, $\Gamma(R)$ does not have isolated vertex. So by Theorem \ref{39}, $\gamma(\Gamma(R))\leq \frac{n-1}{2}$. If $\gamma(\Gamma(R))=\frac{n-1}{2}$, then by Theorem \ref{3k}, $n=3^k$ for $k>1$, which is a contradiction. By Lemma \ref{37}, $\gamma(\Gamma(R))\geq 3$. if $\gamma(\Gamma(R))=3$, then by Lemma \ref{38}, $n=4$, which is impossible. Therefore $3<\gamma(\Gamma(R))<\frac{n-1}{2}$.
$\hfill\Box$
\begin{lem}\label{43} Let $R$ be a non-commutative ring of order $n$ and $Z(R)=\{0\}$. Then $\gamma(\overline{\Gamma}(R))<\frac{n-1}{2}$.
\end{lem}
{\it Proof.}
Since $\overline{\Gamma}(R)$ is a connected graph, $\gamma(\overline{\Gamma}(R))\leq \frac{n-1}{2}$. If $\gamma(\overline{\Gamma}(R))= \frac{n-1}{2}$, then by Theorem \ref{c4}, $\overline{\Gamma}(R)=C_4$ or $\overline{\Gamma}(R)=H\circ K_1$ where $H$ is a connected graph. In both cases $\gamma(\Gamma(R))=2$, contrary to Lemma \ref{37}. Therefore $\gamma(\overline{\Gamma}(R))<\frac{n-1}{2}$.
$\hfill\Box$
\begin{thm} Let $R$ be a non-commutative ring and $Z(R)\neq \{0\}$. Then $\gamma(\overline{\Gamma}(R)) \geq 2$.
\end{thm}
{\it Proof.} On the contrary, let $\gamma(\overline{\Gamma}(R))=1$ and $D=\{x\}$ be a dominating set in $\overline{\Gamma}(R)$. So $C(x)=Z(R) \cup \{x\}$. Let $0\neq y \in Z(R)$. So $x(y+x)=(y+x)x$, and so $y+x \in C(x)$. If $y+x=x$, then $y=0$, which is impossible. Also if $y+x \in Z(R)$, then $x \in Z(R)$, which is a contradiction. Therefore $\gamma(\overline{\Gamma}(R))\geq 2$.
$\hfill\Box$
\begin{cor} If $R$ is a non-commutative ring with unity, then $\Gamma(R)$ has no isolated vertices.
\end{cor}
{\bf Proof of Theorem A.}\\
{\bf i)} Let $\gamma(\Gamma(R)) + \gamma(\overline{\Gamma}(R)) = n $. If $\Gamma(R)$ has no isolated vertex, then $\gamma(\Gamma(R))\leq \frac{n-1}{2}$, by Theorem \ref{39}. Since $\overline{\Gamma}(R)$ is a connected graph, $\gamma(\overline{\Gamma}(R))\leq \frac{n-1}{2}$. Thus $\gamma(\Gamma(R)) + \gamma(\overline{\Gamma}(R))\leq n-1$, which contradicts the assumption. Hence $\Gamma(R)$ has at least one isolated vertex. Thus $\gamma(\overline{\Gamma}(R))=1$ and so $\gamma(\Gamma(R))=n-1$. Therefore $\Gamma(R)=(n-1)K_1$ and $\overline{\Gamma}(R)=K_{n-1}$. By Theorem \ref{18}, $R$ is of type $E$ or $F$. The proof of converse is clear. \\
{\bf ii)} Let $\gamma(\Gamma(R)) + \gamma(\overline{\Gamma}(R)) = n-1$. Since $\overline{\Gamma}(R)$ is a connected graph, $\gamma(\overline{\Gamma}(R))\leq \frac{n-1}{2}$. If $\Gamma(R)$ has no isolated vertex, then by Theorem \ref{39}, $\gamma(\Gamma(R)) \leq \frac{n-1}{2}$. So $\gamma(\Gamma(R)) =\gamma(\overline{\Gamma}(R))= \frac{n-1}{2}$, contrary to Lemma \ref{43}. Hence $\Gamma(R)$ has at least one isolated vertex and so $\gamma(\overline{\Gamma}(R))=1$. Thus $\gamma(\Gamma(R))=n-2$. It is easy to see that $\Gamma(R)=(n-3)K_1\cup P_2$, contrary to Lemma \ref{40}. Therefore $\gamma(\Gamma(R)) + \gamma(\overline{\Gamma}(R))\neq n-1$. \\
{\bf iii)} Let $\gamma(\Gamma(R)) + \gamma(\overline{\Gamma}(R)) = n-2 $. If $\Gamma(R)$ has no isolated vertex, then by Theorem \ref{39}, $\gamma(\Gamma(R)) \leq \frac{n-1}{2}$. Since $\overline{\Gamma}(R)$ is connected graph, $\gamma(\overline{\Gamma}(R)) \leq \frac{n-1}{2}$.
By Lemma \ref{43}, $\gamma(\overline{\Gamma}(R)) < \frac{n-1}{2}$. If $\gamma(\overline{\Gamma}(R)) = \frac{n-1}{2}-1$, then $\gamma(\Gamma(R)) = \frac{n-1}{2}$. By Theorem \ref{3k}, $\Gamma(R) $ is a union of $\frac{n-1}{2}$ copies of $P_2$. It follows that $\overline{\Gamma}(R)$ is a regular multipartite graph of size two. Hence $\gamma(\overline{\Gamma}(R)) =2$, which contradicts the fact that $n$ is not prime. Also if $\gamma(\overline{\Gamma}(R)) \leq \frac{n-5}{2}$, then $\gamma(\Gamma(R)) >\frac{n-1}{2}$, contrary to Theorem \ref{39}. \\
Thus $\Gamma(R)$ has at least one isolated vertex and so $\gamma(\Gamma(R))=n-3$. Let $D=\lbrace x_1,x_2,\ldots, x_{n-3} \rbrace$ be a dominating set in $\Gamma(R)$. Then there exist unique vertices $x_i, x_j\in D$ such that $x_{n-1}\in N(x_i)$ and $x_{n-2}\in N(x_j)$ in $\Gamma(R)$. It is clear that induced subgraph on $D\setminus \{x_i, x_j\}$ is empty. Let $A=\{x_i, x_j, x_{n-1}, x_{n-2} \}$. The proof will be divided into 2 cases.\\
{\bf Case 1.} If $x_i\neq x_j$, then by Lemma \ref{40}, induced subgraph on $A$ is $C_4$, contrary to Lemma \ref{30}. \\
{\bf Case 2.} If $x_i=x_j$, then by Lemma \ref{40}, induced subgraph on $A$ is $K_3$ and so $\Gamma(R)=K_3\cup (n-4)K_1$. Also we has $C(x_{n-1})=\{0, x_{n-1}, x_{n-2}, x_i=x_j\}$. So $4\mid n$. Hence $|R|=n$ is even. The proof of converse is clear.
$\hfill\Box$
\begin{thm} Let R be a non-commutative ring of order $p^2$ and $Z(R)=\{0\}$. Then $\gamma(\Gamma(R))=p+1$.
\end{thm}
{\it Proof.}
We refer the reader to [\cite{A1}, Th. 4].
$\hfill\Box$
\begin{lem}\label{11} Let R be a non-commutative ring of order $p^3$ and $Z(R)=\{0\}$. If $x,y\in V(\Gamma(R))$ and $xy\neq yx$, then $C(x)\cap C(y)=\{0\}$.
\end{lem}
{\it Proof.}
Let $z\in V(\Gamma(R))$. Since $C(z)$ is the addition subgroup of $R$, $|C(z)|\in \{p,p^2\}$. The following 3 cases will be considered.\\
{\bf Case1.} Let $|C(x)|=|C(y)|=p$. Then $|C(x)\cap C(y)|\in \{1,p\}$.\\
If $|C(x)\cap C(y)|=1$, then $C(x)\cap C(y)=\{0\}$.\\
If $|C(x)\cap C(y)|=p$, then $|C(x)\cap C(y)|=|C(x)|$. Since $C(x)\cap C(y)\subseteq C(x)$, $C(x)\cap C(y)=C(x)$. Thus $C(x)\subseteq C(y)$. This contradicts the fact that $y\notin C(x)$.\\
{\bf Case 2.} Let $|C(x)|=|C(y)|=p^2$. Then $|C(x)\cap C(y)|\in \{1,p,p^2\}$.\\
If $|C(x)\cap C(y)|=1$, then $C(x)\cap C(y)=\{0\}$.\\
If $|C(x)\cap C(y)|=p$, then there exist $z\in V(\Gamma(R))$ such that $z\in C(x)\cap C(y)$. By Lemmas \ref{10} and \ref{8}, $C(x)=C(z)$ and $C(y)=C(z)$. Thus $C(x)=C(y)$, a contradiction. \\
If $|C(x)\cap C(y)|=p^2$, then $|C(x)\cap C(y)|=|C(x)|$. Since $C(x)\cap C(y)\subseteq C(x)$, $C(x)\cap C(y)=C(x)$. Thus $C(x)\subseteq C(y)$. This contradicts the fact that $y\notin C(x)$.\\
{\bf Case 3.} Let $|C(x)|=p$ and $|C(y)|=p^2$. Then $|C(x)\cap C(y)|\in \{1,p\}$.\\
If $|C(x)\cap C(y)|=1$, then $C(x)\cap C(y)=\{0\}$.\\
If $|C(x)\cap C(y)|=p$, then $|C(x)\cap C(y)|=|C(x)|$. Since $C(x)\cap C(y)\subseteq C(x)$, $C(x)\cap C(y)=C(x)$. Thus $C(x)\subseteq C(y)$. By Lemma \ref{10}, $C(y)$ is commutative. So for every $z\in C(y)$, $zx=xz$. Thus $z\in C(x)$ and so $C(y)\subseteq C(x)$. Hence $C(y)=C(x)$, which contradicts the fact that $y\notin C(x)$. This completes the proof.
$\hfill\Box$\\
\newline
{\bf Proof of Theorem B.}
Let $x, y \in V(\Gamma(R))$ and $y\notin C(x)$. Then $|C(x)|, |C(y)|$ $\in \{p,p^2\}$. Also by Lemma \ref{11}, $C(x)\cap C(y)=\{0\}$.
Let $|C(x)|=|C(y)|=p$. If $z\in C(x)$, $t\in C(y)$ and $zt=tz$, then by Lemma \ref{8}, $C(x)=C(z)$, $C(y)=C(t)$ and $C(z)=C(t)$. So $C(x)=C(y)$, which is impossible. Therefore $\Gamma(R)$ is the disjoint union of $\ell$ copies of the complete graph of size $p-1$. So $|V(\Gamma(R))|=\ell(p-1)$. On the other hand we have $|V(\Gamma(R))|=p^3-1$. Thus $\ell=p^2+p+1$. Since $\gamma(K_{p-1})=1$, $\gamma(\Gamma(R))=p^2+p+1$, and $(i)$ is proved. \\
Suppose $|C(x)|=p$ and $|C(y)|=p^2$. If $z\in C(x)$, $t\in C(y)$, then $zt\neq tz$. Therefore $\Gamma(R)$ is the disjoint union of $\ell_1$ copies of the complete graph of size $p-1$ and $\ell_2$ copies of the complete graph of size $p^2-1$. So $|V(\Gamma(R))|=\ell_1(p-1)+\ell_2(p^2-1)$. On the other hand we have $|V(\Gamma(R))|=p^3-1$. Thus $\ell_1(p-1)+\ell_2(p^2-1)=p^3-1$ and so $\ell_1+(p+1)\ell_2=p^2+p+1$. Obviously, $\gamma(\Gamma(R))=\ell_1+\ell_2$, and $(ii)$ is proved.
Finally, if $|C(x)|=|C(y)|=p^2$, then $|C(x)+C(y)|=p^4$, which is impossible.
$\hfill\Box$\\
\newline
{\bf Proof of Theorem C.} Let $\gamma(\Gamma(R_i))=m_i$ and $D_i=\{v_{i1}, v_{i2}, \ldots, v_{im_i}\}$ be a dominating set in $\Gamma(R_i)$, where $1\leq i\leq t$. Without loss of generality, let $m_1\leq m_i$, for every $1\leq i\leq t$. We claim that, $D=\{(v_{1j}, 0, \ldots, 0) \mid 1\leq j\leq m_1\}$ is a dominating set in $\Gamma(\prod_{i=1}^{t}R_i)$. For let $(u_1, u_2, \ldots, u_t)\in V(\Gamma(\prod_{i=1}^{t}R_i))$ be an arbitrary vertex. Since $D_1$ is a dominating set in $\Gamma(R_1)$, there exist $v_{1j}\in D_1$ such that $v_{1j}u_1=u_1v_{1j}$ or $v_{1j}=u_1$. Thus $(u_1, u_2, \ldots, u_t)\in N((v_{1j}, 0, \ldots, 0))$ in $\Gamma(\prod_{i=1}^{t}R_i)$. Let $D'=\{(v_{1j}, v_{2j}, \ldots, v_{tj})\mid v_{ij}\in V(\Gamma(R_i)), j=1, \ldots, m_1-1\}\subseteq V(\Gamma(\prod_{i=1}^{t}R_i))$ be a dominating set in $\Gamma(\prod_{i=1}^{t}R_i)$. Let $(u_1, u_2, \ldots, u_t)$ be an arbitrary vertex in $\Gamma(\prod_{i=1}^{t}R_i)$. Since $D'$ is a dominating set in $\Gamma(\prod_{i=1}^{t}R_i)$, there exist a $v_{1j}\in V(\Gamma(R_1))$ such that $u_1v_{1j}=v_{1j}u_1$ or $u_1=v_{1j}$. Whenever $\{v_{11}, v_{12}, \ldots, v_{1(m_1-1)}\}$ is a dominating set in $\Gamma(R_1)$. This contradicts the fact that $D_1$ is a dominating set in $\Gamma(R_1)$. Hence $D$ is a dominating set in $\Gamma(\prod_{i=1}^{t}R_i)$, as claimed. Therefore $\gamma(\Gamma(\prod_{i=1}^{t}R_i))=Min_{1\leq i\leq t}(\gamma(\Gamma(R_i)))$.
$\hfill\Box$
\begin{thm} Let $R_1$ be a non-commutative ring of order $n_1$ and $Z(R_1)=0$. Also, let $R_2$ be a commutative ring of order $n_2$. Then
\begin{center}
$\gamma(\Gamma(R_1\times R_2))=\gamma(\Gamma(R_1))$.
\end{center}
\end{thm}
{\it Proof.}
Let $G=\Gamma(R_1\times R_2)$ and $G'=\Gamma(R_1)\boxtimes K_{n_2}$ such that the members of $V(K_{n_2})$ are the elements of $R_2$. It is easy to see that $G\cong G'$. It is sufficient to prove that $\gamma(G')=\gamma(\Gamma(R_1))$. Let $\gamma(\Gamma(R_1))=m$ and $D=\{v_1, v_2, \ldots, v_m\}$ be a dominating set in $\Gamma(R_1)$. It will be claimed, $D'=\{(v_i, u_0)\mid v_i\in D\}$ is a dominating set in $G'$, where $u_0\in V(K_{n_2})$. Let $(v, u')\in V(G')$ be an arbitrary vertex. Since $D$ is a dominating set in $\Gamma(R_1)$, there exist $v_i\in D$ such that $v_iv=vv_i$ or $v=v_i$. By strongly product structure, $(v, u')\in N((v_i, u_0))$. Let $D''=\{(x_i, u_i)| i=1, \ldots, m-1\}\subseteq V(G')$ be a dominating set in $G'$. Let $(v,u)$ be an arbitrary vertex in $G'$. Since $D''$ is a dominating set in $G'$, there exist a $v_i\in V(\Gamma(R_1))$ such that $(v, u)$ is adjacent to $(v_i, u_i)$. By strongly product structure, $v=v_i$ or $vv_i=v_iv$. whenever $\{v_1, v_2, \ldots, v_{m-1}\}$ is a dominating set in $\Gamma(R_1)$. This contradicts the fact that $D$ is a dominating set in $\Gamma(R_1)$. Thus $D'$ is a dominating set in $G'$, as claimed. Hence $\gamma(G')=\gamma(\Gamma(R_1))$.
$\hfill\Box$
\subsection{Signed domination number in $\Gamma(R)$}
{\bf Proof of Theorem D.} Let $\gamma_s(\Gamma(R))=n-1$. By Lemma \ref{1}, every $v\in V(\Gamma(R))$ is either isolated, an endvertex or adjacent to an endvertex. Let $v\in V(\Gamma(R))$ such that $deg(v)=1$ and $v\in N(u)$. Then $C(v)=\{0, v, u\}$. Thus $u=-v$ and $deg(u)=1$. Hence by Lemma \ref{40}, $\Gamma(R)$ is an empty graph or union of edges.\\
If $n$ is even, then $\Gamma(R)$ containing isolated vertices. Thus $\overline{\Gamma}(R)$ is a complete graph. By Theorem \ref{18}, $R\cong E$ or $R\cong F$ and $(i)$ is proved.\\
If $n$ is odd, then $\Gamma(R)$ is a union of $\frac{n-1}{2}$ copies of $P_2$ and $(ii)$ is proved. The proof of converse is simple.
 $\hfill\Box$
\begin{thm}\label{32} Let $R$ be a non-commutative ring of order $n$, where n is an odd number and $|Z(R)|=c\neq 0$. Then $\gamma_s(\Gamma(R))=n-1$ if and only if $\Gamma(R)$ is the union of $\frac{n-c}{2}$ copies of $P_2$.
\end{thm}
{\it proof.} Let $\gamma_s(\Gamma(R))=n-1$. By Lemma \ref{1}, every $v\in V(\Gamma(R))$ is either isolated, an endvertex or adjacent to an endvertex. If $x$ is an isolated vertex in $\Gamma(R)$, then $C(x)=Z(R)\cup\{x\}$. Thus $O(x)=2$, a contradiction. Let $v\in V(\Gamma(R))$ such that $deg(v)=1$ and $v\in N(u)$. Then $C(v)=Z(R)\cup\{v, u\}$. Since $n$ is an odd number, $u=-v$ and so $deg(u)=1$. Therefore $\Gamma(R)$ is the union of $\frac{n-c}{2}$ copies of $P_2$. The proof of converse is simple.
$\hfill\Box$
\begin{thm} Let the situation be as in Theorem B. Then \\
{\it i)} $\gamma_s(\Gamma(R))=2(p^2+p+1)$.\\
or\\
{\it ii)} $\gamma_s(\Gamma(R))=2(\ell_1+\ell_2)$, where $\ell_1$ and $\ell_2$ satisfy in $\ell_1+(p+1)\ell_2=p^2+p+1$.
\end{thm}
{\it proof.}\\
{\bf i)} As in the proof of Theorem B, $\Gamma(R)$ is the disjoint union of $p^2+p+1$ copies of the complete graphs of size $p-1$. By Theorem \ref{9}, $\gamma_s(K_{p-1})=2$. Hence $\gamma_s(\Gamma(R))=2(p^2+p+1)$, and $(i)$ is proved. \\
{\bf ii)} As in the proof of Theorem B, $\Gamma(R)$ is the disjoint union of $\ell_1$ copies of the complete graphs of size $p-1$ and $\ell_2$ copies of the complete graph of size $p^2-1$, where $\ell_1+(p+1)\ell_2=p^2+p+1$. By Theorem \ref{9}, $\gamma_s(\Gamma(R))=2(\ell_1+\ell_2)$, and $(ii)$ is proved.
$\hfill\Box$
\begin{lem}\label{15} Let R be a non-commutative finite ring and $Z(R)=\{0\}$. If $0\neq |V^-(\overline{\Gamma}(R))|=t$, then the followings are hold.
\item[i)] $\delta(\overline{\Gamma}(R))\leq 2t+1$
\item[ii)] $|R|\leq 4t+2$.
\end{lem}
{\it proof.}\\
 {\bf i)} On the contrary, let $v\in V(\overline{\Gamma}(R))$ and $deg(v)=\delta(\overline{\Gamma}(R))\geq 2t+2$. If $N(v)=\{v_1, \ldots, v_{\delta(\overline{\Gamma}(R))}\}$, then consider the function $f: V(G)\rightarrow \{-1, 1\}$ for which $f(w)=-1$ if and only if $w\in\{v_1, v_2, \ldots, v_{\lfloor\frac{\delta(\overline{\Gamma}(R)}{2} \rfloor}\}$.  Clearly, $f[w]\geq 1$, for all $w\in V(\overline{\Gamma}(R))$. So $f$ is a signed dominating function. This implies that $|V^{-}(\overline{\Gamma}(R))|\geq t+1$, contrary to assumption. \\
{\bf ii)} On the contrary, let $|R|\geq 4t+3$. By Lemma \ref{3}, $\delta(\overline{\Gamma}(R))\geq 2t+2$, which is a contradiction.
$\hfill\Box$
\begin{thm}\label{23} Let $R$ be a non-commutative ring of order $n$ and $Z(R)=\{0\}$. Then $\gamma_s(\overline{\Gamma}(R))\notin \{n-1, n-5\}$.
\end{thm}
{\it proof.} On the contrary, let $\gamma_s(\overline{\Gamma}(R))\in \{n-1, n-5\}$. We consider the following two cases.\\
 {\bf Case 1.} Let $\gamma_s(\overline{\Gamma}(R))=n-1$. If $n=4$, then $\overline{\Gamma}(R)=K_3$ and so $\gamma_s(\overline{\Gamma}(R))\neq n-1$. By corollary \ref{4}, $n\geq 8$. Since $\overline{\Gamma}(R)$ is a connected graph, by Lemma \ref{1}, every $v\in \overline{\Gamma}(R)$ is an endvertex or adjacent to an endvertex. We claim that $\overline{\Gamma}(R)$ have exactly one vertex with degree greater than 1.\\
Suppose that $deg(u), deg(v)>1$. So there are $x,y\in V(\overline{\Gamma}(R))$ such that $x\in N(u)$ and $y\in N(v)$. Since $\overline{\Gamma}(R)$ is a connected graph and by theorem \ref{20}, $diam(\overline{\Gamma}(R))=2$, $v=u$ and $x-v-u$. Hence $\overline{\Gamma}(R)\cong K_{1,n}$, contrary to Theorem \ref{21}. Therefore $\gamma_s(\overline{\Gamma}(R))\neq n-1$.\\
 {\bf Case 2.}  Let $\gamma_s(\overline{\Gamma}(R))=n-5$. Then $|V^{-}(\overline{\Gamma}(R))|=2$. By Lemma \ref{15}, $\delta(\overline{\Gamma}(R)) \leq 5$ and $n\leq 10$. By Corrolary \ref{4}, $n\geq 8$. By Lemmas \ref{26} and \ref{24}, $n\notin\{8, 10\}$.  Let $n=9$ and $v\in V(\overline{\Gamma}(R))$ such that $deg(v)=k=\delta(\overline{\Gamma}(R))$. By Lemma \ref{3}, $k\geq 5$ and so $\delta(\overline{\Gamma}(R))=5$. Hence $|C(v)|\nmid 9$, a contradiction. Therefore $\gamma_s(\overline{\Gamma}(R))\neq n-5$.
$\hfill\Box$
\begin{thm}\label{22} Let $R$ be a non-commutative ring of order $n$ and $Z(R)=\{0\}$. Then $\gamma_s(\overline{\Gamma}(R))=n-3$ if and only if $R$ is isomorphic with one of the following rings\\
$E=\langle x,y~|~2x=2y=0, x^2=x, y^2=y, xy=x, yx=y\rangle$\\
$F=\langle x,y~|~2x=2y=0, x^2=x, y^2=y, xy=y, yx=x\rangle$.
\end{thm}
{\it Proof.}
Let $\gamma_s(\overline{\Gamma}(R))=n-3$. Then $|V^-(\overline{\Gamma}(R))|=1$. By Corrolary \ref{4} and Lemma \ref{15}, $n=4$ and $\delta(\overline{\Gamma}(R))=2$. Thus $\overline{\Gamma}(R)=K_3$. Hence by Theorem \ref{18}, $R$ is one of the following rings\\
$E=\langle x,y~|~2x=2y=0, x^2=x, y^2=y, xy=x, yx=y\rangle$\\
$F=\langle x,y~|~2x=2y=0, x^2=x, y^2=y, xy=y, yx=x\rangle$.\\
The proof of converse is straightforward.
$\hfill\Box$
\begin{cor} Let the situation be as in \ref{22}. If $\gamma_s(\overline{\Gamma}(R))=n-3$, then $\gamma(\overline{\Gamma}(R))=1$.
\end{cor}
{\bf Proof of Theorem E.} \\
Let $v_i\in V(\Gamma(R_i))$ and $deg(v_i)=\delta_i$, where $1\leq i\leq t$. Since $Z(R_i)=\{0\}$, $Z(\prod_{i=1}^{t}R_i)=0$. If $(v_1, v_2, \ldots, v_t)\in V(\Gamma(\prod_{i=1}^{t}R_i))$, then $N[(v_1, v_2, \ldots, v_t)]=\{(x_1, x_2, \ldots, x_t)\mid x_i\in N[v_i]\bigcup\{0\}\}\setminus \{(0, 0, \ldots, 0)\}$. Thus $deg((v_1, v_2, \ldots, v_t))=(\prod_{i=1}^t(\delta_i+2))-2$. Also $(v_1, v_2, \ldots, v_t)$ is a vertex of minimum degree in $\Gamma(\prod_{i=1}^{t}R_i)$. Let $\delta=(\prod_{i=1}^t(\delta_i+2))-2$ and use $u_i\in V(\Gamma(\prod_{i=1}^{t}R_i))$, $1\leq i\leq \delta$, to denote the neighbors of $(v_1, v_2, \ldots, v_t)$. Consider the function $f: V(\Gamma(\prod_{i=1}^{t}R_i))\rightarrow \{-1, 1\}$ for which $f(u_j)=-1$ and $f(u_k)=1$, such that $1\leq j\leq \lfloor\frac{\delta}{2}\rfloor$ and $\lfloor\frac{\delta}{2}\rfloor+1\leq k\leq \delta$. Also, for each $u\in V(\Gamma(\prod_{i=1}^{t}R_i))$ such that $u\neq u_i$, $f(u)=1$. Clearly, $f[w]\geq 1$, for all $w\in V(\Gamma(\prod_{i=1}^{t}R_i))$. Therefore f is a signed dominating function. Hence $|V^{-}(\Gamma(\prod_{i=1}^{t}R_i))|\geq \lfloor\frac{\delta}{2}\rfloor$. Since $\gamma_s(\Gamma(\prod_{i=1}^{t}R_i))=|V(\Gamma(\prod_{i=1}^{t}R_i))|-2|V^{-}(\Gamma(\prod_{i=1}^{t}R_i))|$, it follows that $\gamma_s(\Gamma(\prod_{i=1}^{t}R_i))\leq (\prod_{i=1}^{t}n_i)-1-2\lfloor\frac{\delta}{2}\rfloor$. We consider the following two cases.\\
{\bf Case 1.} Let $\delta$ be odd. Then $\gamma_s(\Gamma(\prod_{i=1}^{t}R_i))\leq (\prod_{i=1}^{t}n_i)-1-2\lfloor\frac{\delta}{2}\rfloor$. It follows that $\gamma_s(\Gamma(\prod_{i=1}^{t}R_i))\leq (\prod_{i=1}^{t}n_i)-1-2\lfloor\frac{\delta-1}{2}\rfloor$. Hence $\gamma_s(\Gamma(\prod_{i=1}^{t}R_i))\leq (\prod_{i=1}^{t}n_i)-\delta$ and $(i)$ is proved.\\
{\bf Case 2.} Let there exist $1\leq i\leq t$ such that $\delta_i$ be even. Then $\gamma_s(\Gamma(\prod_{i=1}^{t}R_i))\leq (\prod_{i=1}^{t}n_i)-1-\delta$ and $(ii)$ is proved.
 $\hfill\Box$


\end{multicols}

\end{document}